\documentstyle{amsppt}
\magnification=1200

\hoffset=-0.5pc
\vsize=57.2truepc
\hsize=38truepc
\nologo
\spaceskip=.5em plus.25em minus.20em
\define\Delt{D}
\define\bv{1}
\define\duality{2}
\define\banach{3}
\define\twilled{4}
\define\koszulon{5}
\define\manipenk{6}
\define\penkoone{7}
\define\Bobb{\Bbb}
\define\do{d^{\nabla}}
\topmatter
\title 
Rinehart complexes and \\
Batalin-Vilkovisky algebras
\endtitle
\author Johannes Huebschmann
\endauthor
\affil 
Universit\'e des Sciences et Technologies de Lille
\\
UFR de Math\'ematiques
\\
F-59 655 VILLENEUVE D'ASCQ C\'edex
\\
Johannes.Huebschmann\@univ-lille1.fr
\endaffil
\endsubjclass
\endkeywords
\endtopmatter
\document
\leftheadtext{Johannes Huebschmann}
\rightheadtext{Rinehart complexes and BV-algebras}

\beginsection Introduction

In this note 
we push further the observation made in \cite\bv\ 
that certain Batalin-Vilkovisky algebras are entirely classical objects:
Let
$(A,L)$ be a Lie-Rinehart algebra such that, as an $A$-module,
$L$ is finitely generated and projective
of finite constant rank.
We will show that the relationship between
generators of the Gerstenhaber bracket 
on $\Lambda_AL$
and connections
on the highest $A$-exterior power of $L$
given in \cite\bv\ 
arises from the canonical pairing between the exterior
$A$-powers of $L$.
Within the framework of duality for Lie-Rinehart algebras
\cite \duality,
this yields another conceptual explanation of the observation,
made already in
\cite\bv\ 
that,
given an exact generator 
for the Gerstenhaber
algebra $\Lambda_AL$,
{\sl the chain
complex underlying the resulting Batalin-Vilkovisky algebra
coincides with the Rinehart complex
computing the corresponding Lie-Rinehart homology\/}.
Thereafter we will spell out some of the connections 
of \cite\bv\ 
with Koszul's paper
\cite\koszulon; in particular we will explain the significance
of the notions of torsionfree linear connection 
and divergence 
in the Lie-Rinehart context.
The reader is assumed familiar with our paper \cite\bv.

\medskip\noindent
{\bf 1. Rinehart complexes and Batalin-Vilkovisky algebras}
\smallskip\noindent
Let $(A,L)$ be a Lie-Rinehart algebra.
By Theorem 1 of \cite\bv,
the formula
$$
a \circ \alpha = a(D\alpha) - \alpha(a),\quad
a \in A,\ \alpha \in L,
\tag1.1
$$
{\sl establishes a bijective correspondence between
right $(A,L)$-connections
$\circ \colon A \otimes_R L \to A$
on $A$ 
(written $(a,\alpha) \mapsto a \circ \alpha$)
and $R$-linear operator $D$ 
on 
$\Lambda_A L$
generating the Gerstenhaber 
bracket\/}.
Under this correspondence,
{\sl flat right $(A,L)$-connections, that is,
right $(A,L)$-module structures on $A$, correspond
to exact operators\/} (i.~e. operators of square zero).
Given 
a right
$(A,L)$-connection
$\circ \colon A \otimes_R L \to A$,
an explicit formula for the generator $D$ is given by 
$$
\aligned
\Delt(\vartheta_1 \wedge \ldots\wedge \vartheta_p) 
&=
\sum_{i=1}^p (-1)^{(i-1)}(1\circ \vartheta_i)
\vartheta_1 \wedge\dots\widehat{\vartheta_i}\ldots\wedge\vartheta_p 
\\
&\quad
+\sum_{1\leq j<k\leq p} (-1)^{(j+k)}
\lbrack\vartheta_j,\vartheta_k \rbrack
\wedge \vartheta_1 \wedge \dots\widehat{\vartheta_j}\dots
\widehat{\vartheta_k}\ldots\wedge \vartheta_p
\endaligned
\tag1.2
$$
where $\vartheta_1,\dots,\vartheta_p \in L$,
cf. (1.6) in \cite\bv.
Furthermore, by 
Theorem 2 of \cite\bv,
given an exact generator 
$\partial$
for the Gerstenhaber
algebra $\Lambda_AL$, 
{\sl the chain
complex underlying the Batalin-Vilkovisky algebra
$(\Lambda_AL,\partial)$
coincides with the Rinehart complex\/}
\linebreak
$(A_{\partial}\otimes_{U(A,L)}K(A,L),d)$
where 
$A_{\partial}$ refers to $A$, endowed with the
right $(A,L)$-module structure determined by $\partial$.
In particular,
the resulting generator $\partial$ 
for the Gerstenhaber bracket
then comes down to a {\sl standard
Koszul-Rinehart operator\/}, and
{\sl when $L$ is projective as an $A$-module,
the Batalin-Vilkovisky algebra
$(\Lambda_AL,\partial)$
computes
the Lie-Rinehart
homology 
$\roman H_*(L,A_\partial) 
$
of $L$ with coefficients in the right $(A,L)$-module\/}
$A_\partial$. 
\smallskip
We now suppose that,
as an $A$-module,
$L$ is  projective  
of finite constant rank (say) $n$,
so that
$\Lambda_A^nL$ is the highest non-zero exterior power of $L$
in the category of $A$-modules.
Write
$
\lambda \colon
L \otimes_R \Lambda_A^nL
@>>>
\Lambda_A^nL
$
for the operation of Lie derivative. For later reference, we recall that,
for 
$\alpha \in L$ and
$x = \xi_1 \wedge \dots \wedge \xi_n \in \Lambda_A^nL$
where $\xi_j \in L$,
$$
\lambda_{\alpha}(x) = \sum_j \xi_1 \wedge \ldots
[\alpha,\xi_j] \ldots \wedge \xi_n .
$$
Theorem 3 in \cite\bv\ 
says that the formula
$$
(1 \circ \alpha) x = \lambda_\alpha (x) - \nabla_\alpha(x)
\tag1.3
$$
establishes
a bijective
correspondence between
right
$(A,L)$-connections
$\circ \colon A \otimes_R L \to A$
on $A$
and $(A,L)$-connections 
$\nabla \colon L \otimes _R \Lambda_A^nL \to \Lambda_A^nL$
on $\Lambda_A^nL$
(written
$(\alpha, x) \mapsto \nabla_\alpha (x) $)
in such a way that
right
$(A,L)$-module structures
on $A$
correspond to left $(A,L)$-module structures on $\Lambda_A^nL$
(i.~e. flat connections).
\smallskip\noindent
{\smc Remark 1.4.}
A conceptual explanation of the formula (1.3) is this:
Given an $A$-module $M$ with an $(A,L)$-connection
$\nabla \colon L \otimes_R M \to M$,
a slight generalization of
\cite\duality\ (2.8) entails that,
given $\phi \in \roman{Hom}_A(\Lambda_A^nL,M)$
and $\alpha \in L$,
the formula
$\phi \circ \alpha = - \lambda^{\nabla}_{\alpha}(\phi)$,
where
$(\lambda^{\nabla}_{\alpha}(\phi))(x) =
\nabla_{\alpha}(\phi(x)) - \phi(\lambda_{\alpha} (x))$
($x \in \Lambda_A^nL$),
endows $\roman{Hom}_A(\Lambda_A^nL,M)$ with
a right $(A,L)$-connection 
$$
\circ \colon 
\roman{Hom}_A(\Lambda_A^nL,M)\otimes_R L 
@>>>
\roman{Hom}_A(\Lambda_A^nL,M).
$$
The formula (1.3)
is a special case thereof,
for $M= \Lambda^n_AL$.
Within the world of supermanifolds,
a version of the association $M \to \roman{Hom}_A(\Lambda_A^nL,M)$
from modules with left connection
to modules with right connection
may be found in \cite\manipenk.
\smallskip
The above observations entail, cf. 
the Corollary in Section 2 of \cite\bv,
that 
the formula
$$
(D\alpha) x = \lambda_\alpha (x) - \nabla_\alpha(x),
\quad x \in \Lambda_A^n L,
\tag1.5
$$
yields
{\sl a bijective correspondence between
$(A,L)$-connections
$\nabla$ on
$\Lambda_A^n L$
and linear operators $D$ generating the Gerstenhaber 
bracket on 
$\Lambda_A L$
in such a way
that flat connections correspond to operators
of square zero, that is, to differentials\/}.
The existence of such a bijective correspondence
is due to
Koszul \cite\koszulon\ 
for the special case where
$A$ is the ring of smooth functions 
and $L$ the $(\Bobb R,A)$-Lie 
algebra of smooth vector fields on a smooth manifold.
The approach in terms of
right $(A,L)$-module structures on $A$
is more general, though,
and yields in particular a conceptual explanation
of the Batalin-Vilkovisky algebra
$(\Lambda_AL,D)$ resulting from a flat connection
$\nabla$ on $\Lambda_A^n L$
as the {\it Rinehart complex\/}
calculating the corresponding {\it Rinehart homology\/}.
It is, furthermore, completely formal and hence immediately generalizes to 
other situations, e.~g. graded Lie-Rinehart algebras,
cf. \cite{\banach,\,\twilled},
vvsheaf versions thereof, 
Lie algebroids on supermanifolds, etc.
\smallskip
We now give another description
of the relationship between
$(A,L)$-connections
$\nabla$ on
$\Lambda_A^n L$
and generators $D$ for the Gerstenhaber bracket
on $\Lambda_AL$.
The canonical pairing
$$
\wedge
\colon
\Lambda_A^* L
\otimes_A
\Lambda_A^{n-*} L
@>>>
\Lambda_A^nL
\tag1.6
$$
of graded $A$-modules
is perfect and its adjoint
$$
\phi\colon\Lambda_A^* L
@>>>
\roman{Hom}_A(\Lambda_A^{n-*} L,\Lambda_A^nL)
=
\roman{Alt}_A^{n-*}(L,\Lambda_A^nL)
\tag1.7
$$
is an isomorphism
of graded $A$-modules. Given $\alpha \in \Lambda_A^* L$,
write $\phi_\alpha \in
\roman{Alt}_A^{n-*}(L,\Lambda_A^nL)$
for the image of $\alpha$
under the isomorphism (1.7).
Thus, when
$\alpha \in \Lambda^p_AL$,
given
$\xi_{p+1},\dots,\xi_n \in L$,
$$
\phi_{\alpha}(\xi_{p+1},\dots,\xi_n)
=\alpha \wedge \xi_{p+1} \wedge\ldots\wedge\xi_n.
$$
For an $(A,L)$-connection
$\nabla \colon M \to \roman{Hom}_A(L,M)$
on a left $A$-module $M$, we denote its
operator of covariant derivative
by
$$
d^{\nabla}
\colon
\roman{Alt}_A(L,M)
@>>>
\roman{Alt}_A(L,M).
$$
The correct degree for an element of
$\roman{Alt}^{n-p}_A(L,\Lambda^n_A L)$ is $p$,
so that (1.7) is degree preserving.
Henceforth we use the standard sign convention for the differential
in a Hom-complex
and, more generally, for 
operators of covariant derivative:
Thus, given a connection
$\nabla \colon L \otimes_R \Lambda^n_A L \to \Lambda^n_AL$,
for $f \in \roman{Alt}^{n-p}_A(L,\Lambda^n_A L)$
and $\xi_p,\dots, \xi_n \in L$,
$$
\alignat1
(\do f)(\xi_p,\dots, \xi_n) 
&=
(-1)^{p+1}
\sum_{p\leq j\leq n}  
(-1)^{j-p}
\nabla_{\xi_j} (f(\xi_p,\dots \widehat \xi_j\dots ,\xi_n))
\\
+
(-1)^{p+1}
&
\sum_{p\leq j<k\leq n} 
(-1)^{j-p+1+k-p+1}
f([\xi_j,\xi_k],
\xi_p,\dots \widehat \xi_j\dots \widehat \xi_k\dots ,\xi_n)
\\
&= 
\sum_{p\leq j\leq n}  
(-1)^{j-1}
\nabla_{\xi_j} (f(\xi_p,\dots \widehat \xi_j\dots ,\xi_n))
\\
&\quad
+
(-1)^{p+1}
\sum_{p\leq j<k\leq n} 
(-1)^{j+k}
f([\xi_j,\xi_k],
\xi_p,\dots \widehat \xi_j\dots \widehat \xi_k\dots ,\xi_n).
\endalignat
$$
The relationship between generators for the Gerstenhaber bracket
and connections on
$\Lambda_A^nL$
is made explicit by the following.

\proclaim{Theorem 1.8}
The relationship
$$
vv\phi_{\Delt(\alpha)} = -
d^{\nabla} (\phi_\alpha),\quad \alpha \in \Lambda^*_A L,
\tag1.8.1
$$
establishes a bijective correspondence
between
generators $\Delt$ for the Gerstenhaber bracket
on
$\Lambda_AL$ and
$(A,L)$-connections $\nabla$ on
$\Lambda_A^nL$
in such a way that
exact generators
correspond to
left $(A,L)$-module structures
$\nabla$,
i.~e. flat $(A,L)$-connections, on
$\Lambda_A^nL$.
This correspondence coincides with the one determined by
{\rm (1.5)}.
\endproclaim

Thus the operator $\Delt$ and the connection $\nabla$
determine each other via (1.8.1).

\demo{Proof}
Let
$\Delt$ be a generator  for the Gerstenhaber bracket
and $\nabla$ a connection  on
$\Lambda_A^nL$.
First we show that
$\Delt$ and $\nabla$ are related by
(1.5) if and only if,
for $0 \leq \ p \leq n$,
the diagram
$$
\CD
\Lambda^p_AL
@>{\phi}>>
\roman{Alt}^{n-p}_A(L,\Lambda^n_A L)
\\
@V{\Delt}VV
@VV{-\do}V
\\
\Lambda^{p-1}_AL
@>{\phi}>>
\roman{Alt}^{n-(p-1)}_A(L,\Lambda^n_A L)
\endCD
\tag1.8.2
$$
is commutative.
In view of what was said above,
given $\xi_p,\dots, \xi_n \in L$,
for any $\alpha \in \Lambda^p_AL$,
$$
\alignat2 
(\do &\phi_{\alpha})
(\xi_p,\dots, \xi_n)
&&=
\sum_{p\leq j\leq n}  
(-1)^{j-1}
\nabla_{\xi_j} (\phi_{\alpha}(\xi_p,\dots \widehat \xi_j\dots ,\xi_n))
\\
&&& 
+
(-1)^{p+1}
\sum_{p\leq j<k\leq n} 
(-1)^{j+k}
\phi_{\alpha}([\xi_j,\xi_k],
\xi_p,\dots \widehat \xi_j\dots \widehat \xi_k\dots ,\xi_n),
\\
&&&=
\sum_{p\leq j\leq n} 
(-1)^{j-1} 
\nabla_{\xi_j}
(\alpha \wedge \xi_p \wedge\dots \widehat \xi_j\ldots\wedge\xi_n)
\\
&&&\quad
+
(-1)^{p+1}
\sum_{p\leq j<k\leq n} 
(-1)^{j+k} 
\alpha \wedge 
[\xi_j,\xi_k]
\wedge\xi_p\wedge\dots\widehat\xi_j\dots\widehat\xi_k\ldots\wedge\xi_n.
\endalignat
$$
In particular, given
$\vartheta_1,\dots ,\vartheta_p \in L$, 
for $\alpha =\vartheta_1 \wedge\ldots\wedge\vartheta_p$,
we obtain
$$
\align
(\do \phi_{\alpha})
&
(\xi_p,\dots, \xi_n)
=
\sum_{p\leq j\leq n} 
(-1)^{j-1}
\nabla_{\xi_j}
(\vartheta_1 \wedge\ldots\wedge\vartheta_p 
\wedge 
\xi_p \wedge\dots \widehat \xi_j\ldots\wedge\xi_n)
\\
&\quad
+
(-1)^{p+1}
\sum_{p\leq j<k\leq n} 
(-1)^{j+k}
\vartheta_1 \wedge\ldots\wedge\vartheta_p \wedge [\xi_j,\xi_k]
\wedge\xi_p\wedge\dots\widehat\xi_j\dots\widehat\xi_k\ldots\wedge\xi_n
\endalign
$$
On the other hand, 
in view of (1.2),
with reference to the 
right $(A,L)$-connection $\circ$
on $A$ which is determined by $D$ and which determines $D$ as well,  
$$
\align
(\phi_{\Delt(\alpha)})
&(\xi_p,\dots, \xi_n)
=(\Delt(\alpha))\wedge\xi_p \wedge\ldots\wedge\xi_n
\\
&=
\sum_{i=1}^p (-1)^{(i-1)}(1\circ \vartheta_i)
\vartheta_1 \wedge\dots\widehat{\vartheta_i}\ldots\wedge\vartheta_p 
\wedge\xi_p \wedge\ldots\wedge\xi_n
\\
&\quad
+\sum_{1\leq j<k\leq p} (-1)^{(j+k)}
\lbrack\vartheta_j,\vartheta_k \rbrack
\wedge \vartheta_1 \wedge \dots\widehat{\vartheta_j}\dots
\widehat{\vartheta_k}\ldots\wedge \vartheta_p
\wedge\xi_p \wedge\ldots\wedge\xi_n
\endalign
$$
Since, as an $A$-module, $L$ is finitely generated
and projective,
it will suffice to 
establish the claim
for the special case where
$
(\vartheta_1, \dots ,\vartheta_p)
=
(\xi_1, \dots ,\xi_p).
$
Indeed, locally $L$ will be free as an $A$-module, and we
may then take 
$\xi_1, \dots ,\xi_n$ to be an $A$-basis.
Now
$$
\alignat2
(\do \phi_{\alpha})
&
(\xi_p,\dots, \xi_n)
&&=
(-1)^{p-1}
\nabla_{\xi_p}
(\xi_1 \wedge\ldots\wedge\xi_p 
\wedge 
\xi_{p+1} \wedge\ldots\wedge\xi_n)
\\
&
\quad
+
(-1)^{p+1}
&&
\sum_{p+1 \leq k \leq n}
(-1)^{p+k} 
\xi_1 \wedge\ldots\wedge\xi_p 
\wedge [\xi_p,\xi_k]
\wedge\xi_{p+1}\wedge\ldots \widehat \xi_k\ldots\wedge\xi_n
\\
&&&=
(-1)^{p+1}
\nabla_{\xi_p}
(\xi_1 \wedge\ldots\wedge\xi_p 
\wedge 
\xi_{p+1} \wedge\ldots\wedge\xi_n)
\\
&
\quad
+
(-1)^{p+1}
&&
\sum_{p+1 \leq k \leq n}
(-1)^k 
[\xi_p,\xi_k]\wedge
\xi_1 \wedge\ldots\wedge\xi_p 
\wedge\xi_{p+1}\wedge\ldots \widehat \xi_k\ldots\wedge\xi_n
\\
(\phi_{\Delt(\alpha)})
&(\xi_p,\dots, \xi_n)
&&
=
(-1)^{(p-1)}(1\circ \xi_p)
\xi_1 \wedge\dots\wedge\xi_{p-1} 
\wedge\xi_p \wedge\ldots\wedge\xi_n
\\
&&& 
+\sum_{1 \leq j \leq p-1}
(-1)^{(j+p)}
\lbrack\xi_j,\xi_p \rbrack
\wedge \xi_1 \wedge \dots\widehat{\xi_j}
\ldots\wedge\xi_{p-1}
\wedge\xi_p \wedge\ldots\wedge\xi_n
\endalignat
$$
Thus 
$$
\alignat2
(-1)^{p-1}(\do \phi_{\alpha})
&
(\xi_p,\dots, \xi_n)
&&=
\nabla_{\xi_p}
(\xi_1 \wedge\ldots\wedge\xi_p 
\wedge 
\xi_{p+1} \wedge\ldots\wedge\xi_n)
\\
&&&
+
\sum_{p+1 \leq k \leq n}
(-1)^k 
[\xi_p,\xi_k]
\wedge
\xi_1 \wedge\ldots\wedge\xi_p 
\wedge\xi_{p+1}\wedge\ldots \widehat \xi_k\ldots\wedge\xi_n
\\
(-1)^{p-1}
(\phi_{\Delt(\alpha)})
&(\xi_p,\dots, \xi_n)
&&
=
(1\circ \xi_p)
\xi_1 \wedge\dots\wedge\xi_{p-1} 
\wedge\xi_p \wedge\ldots\wedge\xi_n
\\
&&& 
+\sum_{1 \leq j \leq p}
(-1)^{j}
\lbrack\xi_p,\xi_j \rbrack
\wedge \xi_1 \wedge \dots\widehat{\xi_j}
\ldots\wedge\xi_{p-1}
\wedge\xi_p \wedge\ldots\wedge\xi_n
\endalignat
$$
where the range of the last summation has been extended to $j \leq p$
(it was $j \leq p-1$ before); 
this extension of summation does not add a non-zero term
since $[\xi_p,\xi_p] = 0$.
\smallskip
In view of (1.3),
for every $\xi \in L$,
with the notation $x =\xi_1 \wedge\ldots\wedge\xi_n$
and $\nabla^{\circ}$
for the $(A,L)$-connection
on $\Lambda^n_AL$
corresponding to $D$ via (1.3),
$$
\align
(1\circ \xi)
\xi_1 \wedge\ldots\wedge\xi_n
&=
\lambda_\xi (x) - \nabla^{\circ}_\xi(x)
\\
&=
\sum_{1 \leq k \leq n}
\xi_1 \wedge\ldots 
\wedge [\xi,\xi_k]\wedge\ldots\wedge\xi_n
-
\nabla^{\circ}_{\xi}(\xi_1 \wedge\ldots\wedge\xi_n) 
\\
&=
\sum_{1 \leq k \leq n}
(-1)^{k-1}[\xi,\xi_k]\wedge
\xi_1 \wedge\ldots \widehat \xi_k\ldots\wedge\xi_n
-
\nabla^{\circ}_{\xi}(\xi_1 \wedge\ldots\wedge\xi_n) .
\endalign
$$ 
Thus we conclude
$$
\alignat2
(-1)^{p-1}
(\phi_{\Delt(\alpha)})
&(\xi_p,\dots, \xi_n)
&&
=
\sum_{1 \leq k \leq n}
(-1)^{k-1}[\xi_p,\xi_k]\wedge
\xi_1 \wedge\ldots \widehat \xi_k\ldots\wedge\xi_n
\\
&&&
\quad
-
\nabla^{\circ}_{\xi_p}(\xi_1 \wedge\ldots\wedge\xi_n)
\\
&&& 
+\sum_{1 \leq j \leq p}
(-1)^{j}
\lbrack\xi_p,\xi_j \rbrack
\wedge \xi_1 \wedge \dots\widehat{\xi_j}
\ldots\wedge\xi_{p-1}
\wedge\xi_p \wedge\ldots\wedge\xi_n
\\
&&&
=
\sum_{p+1 \leq k \leq n}
(-1)^{k-1}[\xi_p,\xi_k]\wedge
\xi_1 \wedge\ldots \widehat \xi_k\ldots\wedge\xi_n
\\
&&&
\quad 
-
\nabla^{\circ}_{\xi_p}(\xi_1 \wedge\ldots\wedge\xi_n)
\\
&&&
=
(-1)^p(d^{\nabla^{\circ}} \phi_{\alpha})(\xi_p,\dots, \xi_n)
\endalignat
$$
Consequently 
the diagram (1.8.2) is commutative
for $0 \leq \ p \leq n$
if and only if
$\nabla = \nabla^{\circ}$, that is,
if and only if
$\Delt$ and $\nabla$ are related by
(1.5).
\smallskip
Reading the calculations backwards we see that, given
an $(A,L)$-connection $\nabla$ 
on $\Lambda^n_AL$,
the operator $D$ determined by the commutativity of the diagram (1.8.2)
is a generator for the Gerstenhaber bracket, and vice versa.
This proves the claim. \qed
\enddemo
\smallskip\noindent
{\smc Remark 1.9.}
With the notation
$\beta = \xi_p \wedge \ldots\wedge \xi_n$,
the Gerstenhaber algebra property entails
$$
0=\Delt(\alpha\wedge\beta)
=
(\Delt\alpha)\wedge\beta
+
(-1)^p
\alpha\wedge(\Delt\beta)
+(-1)^p[\alpha,\beta],
$$
and the reasoning in the above proof may somewhat concisely be summarized
by the formula
$$
0
=
\phi_{\Delt\alpha}\wedge\beta
+
d^{\nabla}\phi_{\alpha}(\beta),
$$
so that
$$
d^{\nabla}\phi_{\alpha}(\beta)
=
(-1)^p
\alpha\wedge(\Delt\beta)
+(-1)^p[\alpha,\beta].
$$

\beginsection 2. The relationship with linear connections

As before, let $(A,L)$ be a Lie-Rinehart algebra, and
maintain the hypothesis that,
as an $A$-module,
$L$ be  projective  
of finite constant rank $n$.
Let $\nabla^L\colon L \otimes_R L \to L$ be an
$(A,L)$-connection on $L$.
Given $\alpha \in L$, define a morphism 
$\Phi^{\nabla^L}_{\alpha}\colon L \to L$
of $A$-modules by means of the formula
$$
\Phi^{\nabla^L}_{\alpha}(\xi)
= [\alpha,\xi] -\nabla^L_{\alpha}(\xi),\quad \xi \in L,
$$
and denote its trace by 
$\roman{Tr}(\Phi^{\nabla^L}_{\alpha}) \in A$.
For $x \in \Lambda^n_A L$, it is obvious that
$$
\roman{Tr}(\Phi^{\nabla^L}_{\alpha})x
=\lambda_\alpha (x) - \nabla_\alpha(x)
=(1 \circ \alpha) x 
=(D\alpha) x;
$$
where
$\nabla\colon L \otimes_R \Lambda^n_A L \to \Lambda^n_A L$ 
is the induced 
$(A,L)$-connection 
on $\Lambda^n_A L$.
Hence,
in view of (1.3) and (1.5) above,
$$
\roman{Tr}(\Phi^{\nabla^L}_{\alpha})x
=(1 \circ \alpha) x 
=(D\alpha) x,
$$
where
\lq\lq $\circ$\rq\rq\ refers to the corresponding
right $(A,L)$-connection on $A$,
and $D$ to the corresponding generator of the Gerstenhaber
bracket on
$\Lambda_A L$.
Consequently,
for any $\alpha \in L$,
$$
\roman{Tr}(\Phi^{\nabla^L}_{\alpha}) 
=1 \circ \alpha =D\alpha .
$$
We note that in \cite\koszulon\ 
the notion of divergence is used which,
given $\alpha \in L$, is written
$
\roman{div}_{\nabla^L}(\alpha)
=\roman{Tr}(\Phi^{\nabla^L}_{\alpha})
$
but we will avoid it since it conflicts with other usages of the terminology
\lq\lq divergence\rq\rq;
see what is said below.
Given the connection $\nabla^L$,
the value
$1 \circ \alpha =D\alpha$
depends only on the
induced connection
$\nabla$ on 
$\Lambda^n_A L$.
The corresponding generator 
$D$ for
the Gerstenhaber
bracket on
$\Lambda_A L$
depends
only on the
induced connection
$\nabla$ on 
$\Lambda^n_A L$
as well.
\smallskip
In \cite\koszulon,
for the special case
where $(A,L)$ is the Lie-Rinehart algebra 
\linebreak
$(C^{\infty}(X), \roman{Vect}(X))$
of
smooth functions and smooth vector fields on a smooth manifold $X$,
Koszul explicitly constructs
a generator $D$ for
the Gerstenhaber
bracket on
$\Lambda_A L$
from a linear connection
on $X$ which has zero torsion.
The zero torsion hypothesis  then entails that,
for any $\alpha$,
$$
\Phi^{\nabla^L}_{\alpha}(\xi)
= [\alpha,\xi] -\nabla^L_{\alpha}(\xi)
=-\nabla^L_{\xi}(\alpha),\quad \xi \in L.
$$
The description (1.5) of the generator $D$ given above
shows that such a generator is determined by
its induced connection 
on $\Lambda^n_A L$
since this description involves only the
induced connection 
on $\Lambda^n_A L$.
The significance of 
the occurrence 
in \cite\koszulon\ 
of a linear connection
on $X$ with zero torsion
is that 
the zero torsion property 
provides an explicit expression
for the corresponding generator 
$\Delt$ entirely within the language of differential operators
which, in turn, is used in \cite\koszulon\ 
to establish the existence of $\Delt$.
The description (1.2) given above which involves a little bit more
homological algebra shows that
the zero torsion hypothesis is superfluous.
This circle of ideas is now closed by the following observation
where we must perhaps make some assumption of the kind which enables us
to apply an appropriate gluing procedure, see below.

\proclaim{Proposition}
Given any 
$(A,L)$-connection $\nabla$ on $\Lambda^n_A L$,
there is an
$(A,L)$-connection on $L$ which induces $\nabla$ and has zero torsion.
\endproclaim

\demo{Proof}
Let $\nabla^L$ be an
$(A,L)$-connection on $L$
having zero torsion,
and let $\nabla$
be the induced
$(A,L)$-connection on $\Lambda^n_A L$.
An arbitrary
$(A,L)$-connection $\overline \nabla$ on $\Lambda^n_A L$
may be written in the form
$$
\overline \nabla_{\alpha}(x) = \nabla_{\alpha} (x) + (\phi(\alpha))x,
\quad
x \in \Lambda^n_A L,
$$
for some  morphism $\phi \colon L \to A$
of $A$-modules.
Let $\Phi \colon L \to \roman{End}_A(L)$
be a  morphism of $A$-modules
such that, for every $\alpha \in L$,
$\phi(\alpha) = \roman{Tr}(\Phi(\alpha))$
and such that
for every $\alpha,\beta \in L$,
$\Phi(\alpha) \beta =\Phi(\beta) \alpha$.
The existence of $\Phi$ is obviouus when
$L$ is a (finitely generated) free $A$-module.
When $L$ is only projective, the existence of $\Phi$ may be
established locally where $L$ is free and then the resulting
data may be glued together appropriately.
For example, $(A,L)$ could be the algebra of smooth functions
on a smooth manifold $W$ and $L$ the space of sections of a
Lie algebroid on $W$.
With this preparation out of the way,
the formula
$$
\overline \nabla^L_{\alpha} = \nabla^L_{\alpha}  + \Phi(\alpha),
\quad
\alpha \in L,
$$
yields   an
$(A,L)$-connection 
$\overline \nabla^L$
on $L$
having zero torsion
which induces
the
$(A,L)$-connection $\overline \nabla$ on $\Lambda^n_A L$.\qed
\enddemo

\noindent
{\smc Remark 1.}
Let $M$ be a free $A$-module
of rank one, with basis element $b$. 
Given
an $R$-module endomorphism 
$\alpha$
of $M$,
its {\it divergence\/} $\roman{div}_b(\alpha) \in A$
is defined by the identity
$$
\alpha (b) =\roman{div}_b(\alpha) b \in M.
$$
In particular,
let $\nabla$ be an
$(A,L)$-connection on $\Lambda^n_AL$
and, given $\alpha \in L$, consider the operation
of generalized Lie derivative
$$
\lambda^{\nabla}_{\alpha}
\colon 
\roman{Alt}^n_A(L,\Lambda^n_AL)
@>>>
\roman{Alt}^n_A(L,\Lambda^n_AL),
\quad
(\lambda^{\nabla}_\alpha\phi)(x) 
=
\nabla_{\alpha}(\phi x) 
-
\phi(\lambda_{\alpha} x),
\quad x \in \Lambda^n_AL,
$$
where
$ \phi \in \roman{Alt}^n_A(L,\Lambda^n_AL)=
\roman{Hom}_A(\Lambda^n_AL,\Lambda^n_AL)$.
Its negative, with $\phi = \roman{Id}$,
yields the right-hand side of (1.3).
Now,
$\roman{Alt}^n_A(L,\Lambda^n_AL)$
is a free $A$-module with basis element
$\phi = \roman{Id}$.
Hence, for $\alpha \in L$,
$$
-\roman{div}_{\phi}(\lambda^{\nabla}_{\alpha})
=\roman{Tr}(\Phi^{\nabla^L}_{\alpha}) 
=1 \circ \alpha =D\alpha .
$$

\noindent
{\smc Remark 2.}
In supermanifold theory, there is also a notion of
{\it integral forms\/}, cf. e.~g. \cite\penkoone.
What corresponds to it under our circumstances
where $(A,L)$ is a Lie-Rinehart algebra
having $L$ finitely generated and projective of finite
constant rank
is the complex $C_L \otimes_UK(A,L) \cong C_L \otimes _A\Lambda_AL$
computing the homology $\roman  H_*(L, C_L)$
where $C_L$ is the dualizing module, cf.
\cite\duality.
When $C_L$ is free as an $A$-module,
this complex of integral forms is just 
a Batalin-Vilkovisky algebra
of the kind
$(\Lambda_AL,\partial)$ described above.
We hope to return to these issues at another occasion.

\bigskip
\widestnumber\key{999}
\centerline{References}
\smallskip\noindent

\ref \no \bv
\by J. Huebschmann
\paper Lie-Rinehart algebras, Gerstenhaber algebras, and Batalin-
Vilkovisky algebras
\jour Annales de l'Institut Fourier
\vol 48
\yr 1998
\pages 425--440
\endref

\ref \no \duality
\by J. Huebschmann
\paper 
Duality for Lie-Rinehart algebras and the modular class
\jour J. reine angew. Math.
\vol 510
\yr 1999
\pages 103--159
\endref

\ref \no \banach
\by J. Huebschmann
\paper Differential Batalin-Vilkovisky algebras arising from
twilled Lie-Rinehart algebras
\paperinfo Poisson Geometry 
\jour Banach Center publications 
\vol 51
\yr 2000
\pages 87--102
\endref

\ref \no \twilled
\by J. Huebschmann
\paper Twilled Lie-Rinehart algebras and differential Batalin-Vilkovisky 
algebras
\paperinfo {\tt math.DG/9811069}
\endref

\ref \no \koszulon
\by J. L. Koszul
\paper Crochet de Schouten-Nijenhuis et cohomologie
\jour Ast\'erisque,
\vol hors-s\'erie,
\yr 1985
\pages 251--271
\paperinfo in E. Cartan et les Math\'ematiciens d'aujourd'hui, 
Lyon, 25--29 Juin, 1984
\endref

\ref \no \manipenk
\by Yu. I. Manin and I. B. Penkov
\paper The formalism of left and right connections on supermanifolds
\paperinfo in: Lectures on supermanifolds, geometrical methods and conformal 
groups, H.-D. Doebner, J. D. Hennig and T. D. Palev, eds.
\pages 3--13
\yr 1989 
\publ World Scientific Publishing Co. Inc.
\publaddr Teaneck, NJ
\endref

\ref \no \penkoone
\by I.B. Penkov
\paper $D$-modules on supermanifolds
\jour Inventiones math.
\vol 71
\yr 1983
\pages 501--512
\endref

\enddocument